\documentclass[12pt]{article}
\usepackage{amsmath}
\usepackage{amssymb}
\usepackage{amsfonts, amsthm}
\begin{document}
\begin{center}{\large\bf Transport in Rotating Fluids
}\end{center} 
\begin{center}{\large\bf Peter Constantin}\end{center}
\begin{center}{Department of Mathematics}\end{center}
\begin{center}{The University of Chicago}\end{center}
\begin{center}{7/23/02}\end{center}
\vspace{1cm}
\noindent{\it Dedicated to Professor M. I. Vishik with deep admiration.}
\newtheorem{thm}{Theorem}
\newtheorem{prop}{Proposition}

\vspace{1cm}
\noindent{\bf Abstract.} {We consider uniformly rotating incompressible Euler and
Navier-Stokes equations. We study the suppression of vertical gradients
of Lagrangian displacement ("vertical" refers to the
direction of the rotation axis). We employ a formalism  that  
relates the total vorticity  to the gradient
of the back-to-labels map (the inverse Lagrangian map, for inviscid
flows, a diffusive analogue for viscous flows). The results include
a nonlinear version of the Taylor-Proudman theorem: in a steady solution
of the rotating Euler equations, two fluid material points which were
initially on a vertical vortex line, will perpetually
maintain their vertical separation unchanged.  For more general situations,
including unsteady flows, we obtain bounds for the vertical gradients of
the Lagrangian displacement that vanish linearly with the maximal local
Rossby number.}

\section{Introduction}

Consider a container filled with water and 
rotated at a constant angular velocity $\Omega$. The 
equations (\cite{chandra}) are
\begin{equation}
\partial_t u + u\cdot\nabla u -\nu\Delta u + \nabla \pi + 2\Omega\widehat{z} \times u = 0
\label{ns}
\end{equation}
with 
\begin{equation}
\nabla\cdot u = 0\label{div}
\end{equation}
$\Delta =\nabla^2$ is the Laplacian.  We use Cartesian coordinates,
the unit vector ${\widehat{z}} = e_3$ is the rotation axis.
The vector 
$u$ is the relative velocity, a function of three space variables and time, 
representing the velocity of the fluid recorded by an ideal observer attached
to the container and rotating thus with it, at uniform angular velocity. The function
$\pi$ contains the physical pressure and the centrifugal force.
The Taylor-Proudman theorem (\cite{chandra}) states that, if one neglects viscosity and inertia, then 
the only time independent solutions are two-dimensional. 

Neglecting viscosity brings us to the three dimensional incompressible 
Euler equations 
\begin{equation}
\partial_t u + u\cdot\nabla u + \nabla\pi + 2\Omega {\widehat{z}}\times u = 0.
\label{re}
\end{equation}
Neglecting inertia brings about the linear equation
\begin{equation} 
\partial_t u + 2\Omega{\widehat{z}}\times u + \nabla\pi  = 0,\quad \nabla\cdot u = 0,
\label{lin}
\end{equation}
and time independent solutions of (\ref{lin}) obey $\partial_3\pi =0$, $u_1 =
-(2\Omega)^{-1}\partial_2 \pi$, $u_2 = (2\Omega)^{-1}\partial_1\pi$. Incompressibility demands then $\partial_3 u_3 = 0$ and standard boundary conditions
require then $u_3 =0$, so indeed, the solutions are two dimensional. 
G.I. Taylor showed experimentally that, when motions are slow and steady,
the motion of the fluid is organized in sheets that remain parallel
to the axis of rotation. Two fluid elements which are initially on a line
parallel to the rotation axis, will remain so. If these elements are a certain
distance apart, they will remain at the same distance apart (\cite{chandra}). 
In Taylor's and in more recent
\cite{sw} experiments the rotation axis coincides with the 
direction of gravity, so we refer to the direction
of the axis as ``vertical''. Thus, for slow, steady inviscid 
motions, there is no vertical transport.  

The main objective of this paper is the quantitative study of the
suppression of vertical transport in the presence of inertia.

Extensive mathematical studies,  (\cite{em}, \cite{bmn} and more recent works),
based on the averaging of the interaction of the fast waves of
(\ref{lin}) with inertia show that two dimensional structures
emerge in the limit $\Omega\to \infty$.  
 
In the present work we employ a formalism  that  
relates the total vorticity  to the gradient
of the back-to-labels map (the inverse Lagrangian map, for inviscid
flows, a diffusive analogue for viscous flows (\cite{c1}, \cite{c2})).
In the presence of rotation the total vorticity decomposes in local
vorticity (curl of relative velocity) and $\Omega e_3$. The ratio between the 
magnitude of the relative  vorticity and $\Omega$ is the local Rossby number
$\rho$. Rotation dominated flows have small Rossby numbers. For steady
solutions of the inviscid equations (\ref{re}), the Lagrangian 
displacement $\lambda (a,t) = X(a,t) -a$ obeys the time independent  
differential equation (see \ref{stdiseq}):
$$
\partial_{a_3}\lambda_3(a,t) + \epsilon_1(a)\partial_{a_3}\lambda_1(a,t)
+ \epsilon_2(a)\partial_{a_3}\lambda_2(a,t) = 0
$$
with 
$$
\epsilon_j(a) = \frac{\rho(a)}{\left (2 + \rho(a)\xi_3(a)\right)}\xi_j(a),
$$
where
$$
\rho (a) = \frac{|\zeta (a)|}{\Omega}, \,\,\, \, \xi(a) = \frac{\zeta(a)}{|\zeta (a)|}
$$
and $\zeta(a)$ is the relative vorticity at $a$. Note that the coefficients
$\epsilon_j(a)$, $j=1, 2$ are small when the local Rossby number 
$\rho(a)$ is small, and vanish when the local direction $\xi(a)$ of the 
relative vorticity is vertical. In these cases we recover the  
statement
$$
\partial_{a_3}\lambda_3 (a,t) = 0,
$$
that is, the Taylor-Proudman theorem.  For more general unsteady flows, we
obtain bounds for vertical gradients of the Lagrangian displacements 
that vanish linearly with $\rho$.

\section{Maps, Velocity and Vorticity}

We consider flow in an open domain  $R$ with smooth boundary $\partial R$ and 
we impose natural boundary conditions: the relative velocity is tangent
to the boundary: $u(x,t)\cdot n(x) = 0$, for $x\in\partial R$. We denote
by $x=(x_1, x_2, x_3)$ or $a=(a_1, a_2, a_3)$ Cartesian coordinates of independent variables. The third Cartesian direction ${\widehat z} = e_3$ is 
singled out. The fluid can be described employing Lagrangian or Eulerian 
descriptions. While doing this, we always keep the Cartesian 
coordinates for independent variables.
The relative Lagrangian path map $X(a,t)$, defined for $a\in R$, and associated to the relative velocity $u(x,t)$ is the solution of the familiar 
ordinary differential equation
\begin{equation}
\partial_t (X(a,t)) = u(X(a,t),t) \label{xeq}
\end{equation}
with initial data
\begin{equation}
X(a,0) = a.\label{xid}
\end{equation}

The rotating Euler equations (\ref{re}) are thus
\begin{equation}
\frac{d}{dt}{(\partial_t X(a,t))} + (\nabla_x\pi)(X(a,t),t)  + 2\Omega{\widehat{z}}\times u(X(a,t),t) = 0\label{elag}
\end{equation}
There is no difference between $\frac{d}{dt}$ and $\partial_t$: the labels
$a$ are held fixed. We use both notations to emphasize the operation at hand.
We take one of the directions $a=(a_1,a_2,a_3)$, say $a_{\alpha}$ and
differentiate the Lagrangian map in the direction $a_{\alpha}$ to obtain the 
vector $\partial_{\alpha}(X(a,t))$. We take the scalar product of the equation (\ref{elag}) with this 
vector, and using the chain rule, we get:
\begin{equation}
\frac{d}{dt}\left (\partial_t X\cdot\partial_{\alpha}X\right ) +
\partial_{\alpha}(\pi -\frac{1}{2}\left|\partial_t X\right |^2) +
2\Omega\left ({\widehat{z}}; \partial_t X; \partial_{\alpha} X\right ) = 0.\label{un}
\end{equation}
where $\left (u;v;w\right)$ denotes the determinant of the matrix whose
columns are $u$, $v$ and respectively,  $w$.
We compute  $\left ({\widehat{z}};\partial _t X ;\partial_{\alpha}X\right) =\partial_t \left ({\widehat{z}}; X ;\partial_{\alpha}X\right) - \left ({\widehat{z}};X ;\partial_{\alpha}\partial_t X\right) = \partial_t \left ({\widehat{z}}; X ;\partial_{\alpha}X\right) -\partial_{\alpha}\left ({\widehat{z}}; X ;\partial_t X\right) + \left ({\widehat{z}};\partial _{\alpha}X ;\partial_{t}X\right)$, and thus, using the antisymmetry of the determinant, we get
\begin{equation}
2\left ({\widehat{z}};\partial _t X ;\partial_{\alpha}X\right)
= \partial_t \left ({\widehat{z}}; X ;\partial_{\alpha}X\right) -\partial_{\alpha} \left ({\widehat{z}}; X ;\partial_{t}X\right).\label{cafe}
\end{equation}
Inserting (\ref{cafe}) in (\ref{un}) we get 
\begin{equation}
\frac{d}{dt}\left\{\left (\partial_t X\cdot\partial_{\alpha}X\right ) + \Omega\left ({\widehat{z}}; X; \partial_{\alpha} X\right ) \right \} +
\partial_{\alpha}\rho  
 = 0\label{central}
\end{equation}
where
\begin{equation}
\rho = \pi - \frac{1}{2}|\partial_t X|^2 - \Omega\left ({\widehat{z}}; X; \partial_{t} X\right ).\label{vienna}
\end{equation}
Integrating in time, from time zero to time $t$, and remembering that time
integration commutes with label derivative, we obtain
$$
\left\{\left (\partial_t X\cdot\partial_{\alpha}X\right ) + \Omega\left ({\widehat{z}}; X; \partial_{\alpha} X\right ) \right \} 
$$
\begin{equation}
= \left\{\left (\partial_t X\cdot\partial_{\alpha}X\right ) + \Omega\left ({\widehat{z}}; X; \partial_{\alpha} X\right ) \right \} _{| t=0} -\partial_{\alpha}\int_0^t\rho .\label{zahlen}
\end{equation}
The initial data area $X(a,0) = 0$ and $\partial_t X (a,0) = u_0(a)$. The
equation (\ref{zahlen}) is therefore 
\begin{equation}
\left \{\partial_t X+ \Omega {\widehat{z}}\times X \right \} \cdot\partial_{\alpha}X
= \left \{u_{0}(a) + \Omega {\widehat{z}}\times a\right \}\cdot e_{\alpha}
-\partial_{\alpha} q\label{musik}
\end{equation}
where 
\begin{equation}
q = \int_0^{t}\rho
\label{applause}
\end{equation}
and $e_{\alpha} = (\delta_{i\alpha})$ is the unit vector associated with the
direction $\alpha$. We consider the ``back-to-labels'' map, (the inverse Lagrangian map) 
\begin{equation}
A(x,t) = X^{-1}(x,t)\label{entracte}
\end{equation}
We read (\ref{applause}) at $a= A(x,t)$, we multiply it by 
$(\partial_i A^{\alpha})(x,t)$ and sum. 
We use the chain rule, and obtain:
$$
u^{i}(x,t) + \Omega \left ({\widehat{z}}; x;e_i\right ) = 
$$
$$
\partial _iA^{\alpha}(x,t)\left \{
u_0^{\alpha}(A(x,t),t) + \Omega \left ({\widehat{z}}; A(x,t); e_{\alpha}\right )\right\} - \partial_ir(x,t),
$$
where 
\begin{equation}
r(x,t) = q(A(x,t),t).\label{asparagus}
\end{equation}
Re-arranging the expressions a bit, we write
$$
u^{i}(x,t) = (\partial_i A^{\alpha}(x,t))u_0^{\alpha}(A(x,t),t) -\partial_i r(x,t)\, +
$$
\begin{equation}
\Omega\left\{({\widehat{z}}; A(x,t), \partial_i A(x,t) - ({\widehat{z}}; x; e_i)\right\}.\label{centragain}
\end{equation}

This is the Weber formula for rotating Euler equations. Obviously, the ``back-to-labels'' map $A$ obeys the equation
\begin{equation}
\partial_t A(x,t) + u(x,t)\cdot\nabla A(x,t) = 0\label{andsalmon}
\end{equation}
which follows by time differentiation from the statement that $A(x,t)$ is the label corresponding to
$x$, $a= A(X(a,t),t)$.
One can verify, by direct calculation, that if $A(x,t)$ solves 
the equation (\ref{andsalmon}) with initial data
\begin{equation}
A(x,0) =x\label{milchrahmstrudel}
\end{equation}
and with velocity computed from (\ref{centragain}), then the velocity
satisfies the rotating Euler equations.
Indeed, writing $\Gamma_0(u,\nabla) = \partial_t + u\cdot\nabla$, and 
applying $\Gamma_0(u, \nabla)$ to (\ref{centragain}) we obtain, using
(\ref{andsalmon})
$$
\Gamma_0(u,\nabla)u^i = -\Omega\left ({\widehat{z}}; u; e_i\right ) -
(\partial_iu^{j})(\partial_j A^{\alpha}(x,t))\left\{u_0^{\alpha}(A) +
\Omega \left ({\widehat{z}}; A(x,t); e_{\alpha}\right )\right \}
$$
\begin{equation}
-\partial_i(\Gamma_0(u,\nabla) r) + (\partial_iu^j)\partial_j r.\label{nocheinkaffe}
\end{equation} 
Now insering the definition (\ref{centragain}) in the equation above, we have
\begin{equation}
\Gamma_0(u,\nabla)u^i = -\Omega\left ({\widehat{z}}; u; e_i\right ) 
-(\partial_i u^{j})u^{j} -\partial_i(\Gamma_0(u,\nabla )r)
- \Omega (\partial_iu^j)\left ({\widehat{z}}; x; e_j\right )\label{bliss}
\end{equation}
We notice that
\begin{equation}
(\partial_iu^j)\left ({\widehat{z}}; x; e_j\right ) = 
\partial _i\left ({\widehat{z}}; x; u\right) - \left ({\widehat{z}}; e_i; u\right)\label{spazieren}
\end{equation}
so that the equation (\ref{bliss}) reads
\begin{equation}
\Gamma_0(u,\nabla) u^i = -2\Omega\left ({\widehat{z}}; u, e_i\right ) -\partial_i \pi\label{zahlenbitte}
\end{equation}
with
\begin{equation}
\pi = \Gamma_0(u,\nabla) r + \frac{1}{2}|u|^2 + \Omega \left ({\widehat{z}}; x; u\right )\label{dankeschon}
\end{equation}
The obtained equation (\ref{zahlenbitte}) is just (\ref{re}), so we have proved therefore:
\begin{thm} Any twice continuously differentiable solution of the active 
vector system
(\ref{centragain}, \ref{andsalmon}) solves the rotating Euler equations
(\ref{re}) with pressure $\pi$ given in (\ref{dankeschon}).
Viceversa, if $u$ is a twice differentiable solution of the rotating 
Euler equations, then it obeys
(\ref{centragain}) with $A$ determined by (\ref{andsalmon}).
\end{thm}
The Navier-Stokes equations (\ref{ns}) in a rotating frame admit a somewhat
similar treatment. 
One considers a map $A(x,t)$ that obeys
\begin{equation}
\left (\partial_t + u \cdot\nabla -\nu\Delta\right )A(x,t) = 0
\label{nua}
\end{equation}
together with the initial condition (\ref{milchrahmstrudel}). It is convenient
to denote 
\begin{equation}
\Gamma_{\nu}(u,\nabla) = \partial_t + u \cdot\nabla -\nu\Delta\label{gammanu}
\end{equation}
the operator of advection by $u$ and diffusion with diffusivity $\nu$.
One starts with an ansatz like (\ref{centragain}),
\begin{equation}
u^i = (\partial_iA^{\alpha})v^{\alpha} + \Omega\left\{({\widehat {z}};A;\partial_i A) - ({\widehat{z}}; x, e_i)\right\} - \partial_ir \label{webns}
\end{equation}
which is exactly (\ref{centragain}) except that instaed of $u_{0}(A(x,t))$
we have now $v(x,t)$. Note that, in the inviscid case,
the function $v_0(x,t)=u_0(A(x,t))$ obeys $\Gamma_0(u,\nabla)v_0 =0$.
The viscosity introduces corrections to this equation. We apply
$\Gamma_{\nu}(u, \nabla)$ to (\ref{webns}) and compute the various pieces
separately. First, we get
$$
\Gamma_{\nu}(u,\nabla)\left\{\left (\partial_iA^{\alpha}\right )v^{\alpha}\right\} = 
$$
\begin{equation}
(\partial_i A^{\alpha})\Gamma_{\nu}(u,\nabla )v^{\alpha} -2\nu(\partial_j\partial_i A^{\alpha})\partial_jv^{\alpha} - (\partial_i u^j)\left (\partial_jA^{\alpha}\right )v^{\alpha}.\label{vone}
\end{equation}
Secondly, we have
\begin{equation}
\Gamma_{\nu}(u,\nabla)(\partial_i r) = \partial_i \left (\Gamma_{\nu}(u,\nabla)r\right ) - (\partial_iu^j)(\partial_j r)\label{vtwo}
\end{equation}
Subtracting (\ref{vtwo}) from (\ref{vone}) we get
$$
\Gamma_{\nu}(u,\nabla)\left\{\left (\partial_iA^{\alpha}\right )v^{\alpha} - \partial_i r\right\} = -\partial_i\left (\Gamma_{\nu}(u,\nabla) r\right ) +
$$
\begin{equation}
(\partial_i A^{\alpha})\Gamma_{\nu}(u,\nabla )v^{\alpha} -2\nu(\partial_j\partial_i A^{\alpha})\partial_jv^{\alpha} - (\partial_i u^j)\left\{\left (\partial_jA^{\alpha}\right )v^{\alpha} -\partial_j r\right\}.\label{vthree}
\end{equation}
Let us give a temporary name to the term involving $\Omega$ in (\ref{webns}):
\begin{equation}
U^i= \Omega\left \{({\widehat{z}}; A; \partial_i A) - ({\widehat{z}}; x; e_i)\right\}. \label{tempu}
\end{equation}
Then, using (\ref{webns}) we write (\ref{vthree}) as
$$
\Gamma_{\nu}(u,\nabla )u = (\partial_i A^{\alpha})\Gamma_{\nu}(u,\nabla )v^{\alpha} -2\nu(\partial_j\partial_i A^{\alpha})\partial_jv^{\alpha} -\partial_i\left (\Gamma_{\nu}(u,\nabla) r\right )
$$
\begin{equation}
- (\partial_iu^j)u^j + (\partial_iu^j)U^j + \Gamma_{\nu}(u,\nabla )U^i
\label{vfour}
\end{equation}
Now we compute $\Gamma_{\nu}(u,\nabla)U^i$:
\begin{equation}
\Gamma_{\nu}(u,\nabla )U^i = -\Omega\left \{(\partial_i u^j)({\widehat{z}}; A; \partial_j A) + ({\widehat{z}}; u; e_i) + 2\nu ({\widehat{z}}; \partial_j A,\partial_j\partial_i A)\right\}\label{gamaui}
\end{equation}
Using
$$
\partial_i\left (-\Omega ({\widehat{z}}; x; u)\right ) = -\Omega({\widehat{z}};x; (\partial_i u^j)e_j) - \Omega ({\widehat{z}}; e_i; u)
$$
and noticing that two terms $ \Omega (\partial_i u^j)({\widehat{z}}; A; \partial_j A)$ cancel, we deduce that
\begin{equation}
(\partial_iu^j) U^j + \Gamma_{\nu}(u,\nabla )U^i = 
-\Omega \partial_i({\widehat{z}}; x; u) - 2\Omega({\widehat{z}}; u; e_i)
-2\Omega\nu({\widehat{z}};\partial_j A; \partial_i\partial_j A).\label{uiuj}
\end{equation}
Therefore, we have arrived at
$$
\Gamma_{\nu}(u, \nabla)u = -\partial_i\left\{\Gamma_{\nu}(u,\nabla) r + \frac{1}{2}|u|^2 + \Omega({\widehat{z}}; x; u)\right\} - 2\Omega({\widehat{z}}; u; e_i)  +
$$
\begin{equation}
(\partial_iA^{\alpha})\Gamma_{\nu}(u,\nabla )v^{\alpha} -2\nu(\partial_j\partial_iA^{\alpha})(\partial_j v^{\alpha}) - 2\Omega\nu({\widehat{z}};\partial_j A; \partial_i\partial_j A).\label{vfive}
\end{equation}
Now the last piece in (\ref{vfive}) vanishes if
$v$ solves the equation
\begin{equation}
\Gamma_{\nu}(u,\nabla v) = 2\nu C^{\alpha}_{j;\beta}v^{\alpha}_j +
2\Omega\nu ({\widehat{z}}; \partial_j A; C^{.}_{j;\beta})\label{veq}
\end{equation}
Here
\begin{equation}
C^{\alpha}_{j;\beta} = (\nabla A)^{-1}_{k \beta}\partial_j\partial_k A^{\alpha}
= - \Gamma^{\alpha}_{\beta\gamma}(\partial_jA^{\gamma})\label{cdef}
\end{equation}
is related to the Riemann-Christoffel symbol $\Gamma^{\alpha}_{\beta \gamma}.$
Thus we have proved
\begin{thm} Assume that $A$ solves (\ref{nua}) with $u$ given by
(\ref{webns}) and with $v$ evolving according to (\ref{veq}).
Then $u$ solves the rotating Navier-Stokes equation (\ref{ns}).
Vice-versa, if $u$ solves the rotating Navier-Stokes equations
and if we solve the linear equation (\ref{nua}) and then the
linear inhomogeneous equation (\ref{veq}), then (\ref{webns}) holds.
\end{thm}
Let us derive the analogue of (\ref{webns}) for the vorticity. We
will take the curl of (\ref{webns}) and express the vorticity
$\omega = \nabla\times u$ in terms of the diffusive back-to-labels
map $A$ associated to the relative velocity, and the virtual vorticity
$\zeta = \nabla_A\times v$. This object is a generalization 
of the Cauchy invariant $\zeta_0(x,t) = \omega_0(A(x,t))$, to which it 
reduces if $\nu = 0$. Thus, $\zeta$ is defined by
\begin{equation}
\zeta^{\alpha}(x,t) = \frac{1}{2}\epsilon_{\alpha\beta\gamma}(\nabla A)^{-1}_{j\beta}(\partial_j v^{\gamma}).\label{zeta}
\end{equation}
We compute the curl of (\ref{webns}); we start with the curl of $(\partial_iA^{\alpha})v^{\alpha}$:
$$
\epsilon_{ijk}\partial_j\left\{(\partial_kA^{\alpha})v^{\alpha}\right\} =
 \epsilon_{ijk}(\partial_k A^{\alpha})\partial_jv^{\alpha}.
$$
Using the convenient notation $v^{\alpha}_{;\beta} = (\nabla A)^{-1}_{j\beta}\partial_jv^{\alpha}$, we obtain
$$
\epsilon_{ijk}\partial_j\left\{(\partial_kA^{\alpha})v^{\alpha}\right\} =
\epsilon_{ijk}(\partial_k A^{\alpha})(\partial_jA^{\beta})v^{\alpha}_{;\beta} = \frac{1}{2}\epsilon_{ijk}(\partial_k A^{\alpha})(\partial_jA^{\beta})(v^{\alpha}_{;\beta} - v^{\beta}_{;\alpha})
$$
The last equality holds because the contribution of 
$\frac{1}{2}(v^{\alpha}_{;\beta} + v^{\beta}_{;\alpha})$ vanishes. Thus
$$
\epsilon_{ijk}\partial_j\left\{(\partial_kA^{\alpha})v^{\alpha}\right\} =
\frac{1}{2}\epsilon_{ijk}(\partial_kA^{\alpha})(\partial_jA^{\beta})\epsilon_{\beta\alpha\gamma}\zeta^{\gamma}
$$
and consequently
\begin{equation}
\epsilon_{ijk}\partial_j\left\{(\partial_kA^{\alpha})v^{\alpha}\right\} = 
\frac{1}{2}\epsilon_{ijk}(\partial_j A; \partial_k A; \zeta)\label{cone}
\end{equation}
Now we compute the curl of (\ref{tempu}):
$$
\epsilon_{ijk}\partial_jU^k = \Omega\epsilon_{ijk}\left\{({\widehat{z}};\partial_j A;\partial_k A;) + ({\widehat{z}}; A; \partial_j\partial_k A) - ({\widehat{z}}; e_j; e_k)\right\}
$$
which gives
\begin{equation} 
\epsilon_{ijk}\partial_jU^k = \Omega\epsilon_{ijk}\left\{({\widehat{z}};\partial_j A;\partial_k A;) - ({\widehat{z}}; e_j; e_k)\right\}.\label{ctwo}
\end{equation}
Adding (\ref{cone}) to (\ref{ctwo} we deduce
\begin{equation}
\omega^{i} = \frac{1}{2}\epsilon_{ijk}\left\{(\partial_jA;\partial_k A; \zeta)
+ 2\Omega\left [(\partial_j A;\partial_k A; {\widehat{z}}) - (e_j; e_k; {\widehat{z}})\right]\right\}.\label{cauchyrnu}
\end{equation}
This states that the total vorticity $\omega + 2\Omega{\widehat{z}}$ obeys
the Cauchy invariance (\cite{c1}, \cite{c2})
\begin{equation}
\omega + 2\Omega{\widehat{z}} = {\mathcal C}(\nabla A; \zeta + 2\Omega{\widehat{z}})\label{cauchytot}
\end{equation}
with respect to the trajectories of the relative velocity. Here
we used the notation ${\mathcal C}(M; v)= \frac{1}{2}\epsilon_{ijk}(M_{.j}; M_{.k}; v)$.

\section{Lagrangian Transport} 
We will look more closely now at the term involving $\Omega$ in the expression
(\ref{cauchyrnu}). We start by writing
\begin{equation}
\omega = {\mathcal C}(\nabla A; \zeta) + 2\Omega {\mathcal R}(\ell),\label{r}
\end{equation}
Here
\begin{equation}
\ell (x,t) = A(x,t) -x\label{l}
\end{equation}
is the Eulerian-Lagrangian displacement, and
\begin{equation}
\left ({\mathcal R}(\ell)\right )_i  = \frac{1}{2}\epsilon_{ijk}\left\{({\widehat{z}};\partial_j \ell; \partial_k\ell) + ({\widehat{z}}; e_j; \partial_k\ell) + 
({\widehat{z}}; \partial_j\ell; e_k)\right\}\label{rone}
\end{equation}

Using ${\widehat{z}} = e_3$, the term depending linearly on $\ell$ is
$$
\frac{1}{2}\epsilon_{ijk}\left\{(e_j; \partial_k \ell; e_3) + (\partial_j \ell; e_k; e_3)\right\} = \epsilon_{ijk}(e_j; \partial_k\ell; e_3) =
$$
$$
\{-\delta_{k1}\delta_{j2}\delta_{i3}(\partial_2 \ell; e_1;e_3)
+ \delta_{k1}\delta_{j3}\delta_{i2}(\partial_3\ell; e_1; e_3) +
\delta_{k2}\delta_{j1}\delta_{i3}(\partial_1\ell; e_2; e_3) 
$$
$$
- \delta_{k2}\delta_{j3}\delta_{i1}(\partial_3; e_2; e_3) \} 
= \left \{(\partial_2\ell_2) e_3 -(\partial_3\ell_2)e_2 \right\} +
$$
\begin{equation}
\left\{ (\partial_1\ell_1)e_3 - (\partial_3\ell_1)e_1\right\} =   \left\{ -\partial_z\ell + (\nabla\cdot\ell)e_3\right\}
\label{lone}
\end{equation}
The quadratic term is
$$
\frac{1}{2}\epsilon_{ijk}(\partial_j\ell; \partial_k\ell; e_3) =
\frac{1}{2}\epsilon_{ijk}\epsilon_{\alpha\beta\gamma}(\partial_j\ell_{\alpha})(\partial_k\ell_{\beta})\delta_{\gamma 3} = 
$$
\begin{equation}
= \epsilon_{ijk}(\partial_j\ell_1)(\partial_k\ell_2) = \nabla\ell_1\times\nabla\ell_2.\label{ltwo}
\end{equation}
Putting together (\ref{lone}) and (\ref{ltwo}) we have thus
\begin{equation}
{\mathcal R}(\ell) = -\partial_z\ell + (\nabla\cdot\ell){\widehat{z}} + (\nabla\ell_1\times\nabla\ell_2)\label{rl}
\end{equation}
The detailed expression is:

\begin{equation}
{\mathcal R}(\ell) = \left (
\begin{array}{c} -\partial_3\ell_1 + (\partial_2\ell_1)(\partial_3\ell_2) -
 (\partial_3\ell_1)(\partial_2\ell_2) \\ -\partial_3\ell_2 + (\partial_3\ell_1)(\partial_1\ell_2) - (\partial_1\ell_1)(\partial_3\ell_2) \\ -\partial_3\ell_3
 + (\nabla\cdot \ell) + (\partial_1\ell_1)(\partial_2\ell_2) - (\partial_2\ell_1)(\partial_1\ell_2)
\end{array}
\right )
\label{rlone}
\end{equation}

Let us consider the case of the rotating three dimensional Euler equations. In this case the back-to-labels map is the inverse Lagrangian
path map. The Eulerian-Lagrangian displacement $\ell(x,t)$ is related to
the Lagrangian displacement
\begin{equation}
\lambda (a,t) = X(a,t) - a\label{lagdis}
\end{equation}
by
\begin{equation}
\lambda (a,t) + \ell(x,t) = 0, \,\, \, \, a = A(x,t). \label{lal}
\end{equation}
As long as the evolution is smooth, the
incompressibility constraint ($\nabla\cdot u = 0$) implies
\begin{equation}
Det(\nabla A) = 1 \label{det3}
\end{equation}
In view of $\nabla A = {\mathbf {1}} + (\nabla\ell)$, one has therefore
$$
\nabla\cdot\ell + (\partial_2\ell_2)(\partial_3\ell_3) -(\partial_3\ell_2)(\partial_2\ell_3) + (\partial_1\ell_1)(\partial_3\ell_3) - (\partial_3\ell_1)(\partial_1\ell_3) \, + 
$$
\begin{equation}
(\partial_1\ell_1)(\partial_2\ell_2) - (\partial_2\ell_1)(\partial_1\ell_2) = 0
\label{threed}
\end{equation}
Now comes the main observation concerning ${\mathcal R}(\ell)$: The only terms
in (\ref{rlone}) that are not explicitly multiples of some component of
$\partial_z\ell$ are found in the third component; but, in view of (\ref{threed}),
$(\nabla\cdot\ell) + (\partial_1\ell_1)(\partial_2\ell_2) -(\partial_2\ell_1)(\partial_1\ell_2)$ is a quadratic expression involving exclusively
multiples of $\partial_z\ell$. Therefore one can factor
\begin{equation}
{\mathcal R}(\ell) = -{\mathcal M}(\nabla \ell)\partial_z\ell, \label{mrl}
\end{equation}
where the matrix ${\mathcal M}(\nabla\ell)$ is:
\begin{equation}
{\mathcal M}(\nabla\ell) = \left ( 
\begin{array}{ccc}
1 + \partial_2\ell_2 & - \partial_2\ell_1 & 0 \\
-\partial_1\ell_2 & 1 + \partial_1\ell_1 & 0 \\
- \partial_1\ell_3 & - \partial_2\ell_3  & 1 + \partial_1\ell_1 + \partial_2\ell_2
\end{array}
\right )\label{mel}
\end{equation}
Returning to (\ref{r}) we obtain the relation
\begin{equation}
{\mathcal M}(\nabla \ell)\partial_z\ell = \frac{1}{2\Omega}\left \{{\mathcal C}(\nabla A; \zeta) - \omega \right\}\label{mdz}
\end{equation}
Let us consider the matrix
\begin{equation}
{\mathcal N}(\nabla\ell) = \left (
\begin{array}{ccc}
1 + \partial_1\ell_1 & \partial_2\ell_1 & 0 \\
\partial_1\ell_2 & 1+ \partial_2\ell_2 & 0 \\
0 & 0 & 1
\end{array}
\right )
\label{nel}
\end{equation}
and multiply both sides of (\ref{mdz}) by (\ref{nel}) from the left. 
Let us denote the right hand side by the single letter $s$ (for small).
\begin{equation}
s = \frac{1}{2\Omega}{\mathcal N}(\nabla\ell)\left \{{\mathcal C}(\nabla A; \zeta) - \omega \right\}
\label{s}
\end{equation}
We obtain thus
\begin{equation}
D_2(\nabla\ell)\partial_z\ell_1 = s_1, \label{s1}
\end{equation}
\begin{equation}
D_2(\nabla\ell)\partial_z\ell_2 = s_2, \label{s2}
\end{equation}
and
\begin{equation}
(1+ t_2(\nabla\ell))\partial_z\ell_3 - (\partial_1\ell_3)(\partial_3\ell_1)
- (\partial_2\ell_3)(\partial_3\ell_2) = s_3. \label{s3}
\end{equation}
We use the notations
\begin{equation}
D_2(\nabla\ell) = (1 + \partial_1\ell_1)(1 + \partial_2\ell_2) - (\partial_1\ell_2)(\partial_2\ell_1),\label{D2}
\end{equation}
\begin{equation}
d_2 (\nabla\ell) =  (\partial_1\ell_1)(\partial_2\ell_2) - (\partial_1\ell_2)(\partial_2\ell_1),\label{d2}
\end{equation}
\begin{equation}
t_2(\nabla\ell) = \partial_1\ell_1 + \partial_2\ell_2.
\label{t2}
\end{equation}
From definitions  obviously
\begin{equation}
D_2(\nabla\ell) = 1 + t_2(\nabla\ell)  + d_2(\nabla\ell). \label{eqobv}
\end{equation}
Also, the information contained in (\ref{threed}) can be written as
\begin{equation}
t_2(\nabla\ell) + d_2(\nabla\ell) + (1+ t_2(\nabla\ell))\partial_z\ell_3
= (\partial_1\ell_3)(\partial_3\ell_1) + (\partial_2\ell_3)(\partial_3\ell_2)
\label{threede}
\end{equation}
Substituting (\ref{threede}) in (\ref{s3}) we obtain
\begin{equation}
t_2(\nabla\ell) + d_2(\nabla\ell)  = - s_3\label{t2+d2}
\end{equation}
Using (\ref{eqobv}) we obtain 
\begin{prop}
Let $\omega $ be the relative vorticity $\omega = \nabla\times u$ 
of a solution of the rotating three dimensional Euler equations (\ref{re}), 
let $A$ be the inverse Lagrangian map associated to the relative velocity 
$u$, let $\zeta$ be the Cauchy invariant ($\zeta(x,t) = \omega_0(A(x,t))$
, with $\omega_0 = \omega _{|t=0}$), and let $D_2(\nabla\ell)$ be the determinant of the block $(1,2)$ of $\nabla A$,
i.e.
\begin{equation}
D_2(\nabla\ell) = (\partial_1 A_1)(\partial_2A_2) - (\partial_1A_2)(\partial_2A_1).\label{D2a}
\end{equation}
Then one has
\begin{equation}
D_2(\nabla\ell) = 1 + \frac{1}{2\Omega}\left \{\omega _3 - \left (\partial_1 A; \partial_2 A; \zeta\right )\right \},\label{lowb}
\end{equation}
and the equations (\ref{s1}, \ref{s2}, \ref{s3}).
\end{prop}
Let us assume now that 
\begin{equation}
\sup_{t\in [0,T]}\|\omega (\cdot, t)\|_{L^{\infty}(dx)} = M <\infty .
\label{m}
\end{equation}
Let us assume also that on a time interval $I$
we have 
\begin{equation}
\sup_{t\in I}\|\nabla\ell (\cdot, t)\|_{L^{\infty}(dx)}\le g.\label{g}
\end{equation}
Clearly then
\begin{equation}
|s_3| \le (1+ 2g + 3g^2)\frac{M}{\Omega}\label{s_3}
\end{equation}
and
\begin{equation}
|s_j| \le (1+2g)(1+ 2g + 3g^2)\frac{M}{\Omega}\label{s_j}
\end{equation}
$j = 1, 2$ hold on this time interval. We deduce from (\ref{lowb}) that
\begin{equation}
D_2(\nabla \ell) \ge 1 - (1+ 2g + 3g^2)\frac{M}{\Omega}\label{lbd2}
\end{equation}
holds on the time interval $I$. 
Using (\ref{s1}, \ref{s2}) we deduce that
\begin{equation}
|\partial_z\ell_j| \le (1+ 2g)(1 + 2g + 3g^2)\left (1 - (1+2g+3g^2)\frac{M}{\Omega}\right )^{-1}\frac{M}{\Omega}\label{dzlj}
\end{equation}
holds for $j = 1, 2$ on the time interval $I_i$. 
From (\ref{eqobv}, \ref{lbd2}) we deduce that
\begin{equation}
1 + t_2(\nabla\ell) \ge  1 - 2g^2 - (1+ 2g + 3g^2)\frac{M}{\Omega}\label{lb3}
\end{equation}
holds on the same interval. From (\ref{s3}, \ref{lb3}) and from (\ref{dzlj})
we deduce now that
\begin{equation}
|\partial_z\ell_3| \le C_g \frac{M}{\Omega}
\label{dzl3}
\end{equation}
with
$$  
C_g = \left (1 - 2g^2 - (1+ 2g + 3g^2)\frac{M}{\Omega}\right )^{-1}\times
$$
\begin{equation}
\left [1 + 2g (1+2g)\left (1 - (1+2g+3g^2)\frac{M}{\Omega}\right )^{-1}\right ](1+ 2g + 3g^2)\frac{M}{\Omega}\label{cg}
\end{equation}
In order to simplify the exposition, let us assume for example, that
\begin{equation}
g \le \frac{1}{4}\label{gspec}
\end{equation}
Let us introduce the maximal local Rossby number
\begin{equation}
\rho = \frac{M}{\Omega} = \frac{\sup_{0\le t\le T}\|\omega( \cdot,t)\|_{L^{\infty}(dx)}}{\Omega}.
\label{ros}
\end{equation}
The most stringent condition of invertibility, $1>2g^2 + (1+2g + 3g^2)\rho$
is satisfied if the Rossby number is small enough.
Let us assume that the Rossby number satisfies
\begin{equation}
\rho \le \frac{1}{4}\label{rospec}
\end{equation}
Then we deduce
\begin{prop}
Let a solution of the rotating Euler equation (\ref{re}) satisfy  
\begin{equation}
\sup_{t\in I} \|\nabla\ell (\cdot, t)\|_{L^{\infty}}\le \frac{1}{4} 
\label{gspecs}
\end{equation}
on a time interval $I$ and assume that the local Rossby number $\rho$
(\ref{ros}) does not exceed $1/4$ (\ref{rospec}). 
Then
\begin{equation}
\sup_{t\in I}\|\partial_z\ell\|_{L^{\infty}(dx)}\le 14\rho
\label{dzlbound}
\end{equation}
holds on the interval $I$.
\end{prop}
In order to bound the vertical derivative of the direct Lagrangian map
$X(a,t)$ we recall that
\begin{equation}
(\nabla_a X)(A(x,t),t) = \left (\nabla A (x,t)\right )^{-1}\label{nabalx}
\end{equation}
and in particular
\begin{equation}
\frac{\partial X}{\partial a_3}(A(x,t),t) = (\nabla A_1(x,t))\times (\nabla A_2(x,t)).\label{nabl}
\end{equation}
Expressing this in terms of the displacement $\ell$ we arrive at
\begin{equation}
 \frac{\partial X}{\partial a_3}(A(x,t),t) =  e_3 + \left (
\begin{array}{c}
(\partial_2 \ell_1)(\partial_3\ell_2) - 
(1 + \partial_2\ell_2)(\partial_3\ell_1) \\ 
(\partial_1\ell_2)(\partial_3\ell_1) - 
(1+ \partial_1\ell_1)(\partial_3\ell_2)\\
t_2(\nabla\ell) + d_2(\nabla \ell)
\end{array}
\right )\label{dx}
\end{equation}
In view of (\ref{eqobv}, \ref{t2+d2}, \ref{lowb}, \ref{s_3}, \ref{dzlj}) 
we deduce
\begin{prop}
Let a solution of the rotating Euler equation (\ref{re}). Consider the
Lagrangian map $X(a,t)$ associated to the relative velocity.
Then
\begin{equation}
\frac{\partial X_3}{\partial a_3}(a,t) = 1 + \frac{1}{2\Omega}\left (\omega_3(X(a,t),t) - \zeta (a)\cdot \frac{\partial X}{\partial a_3}(a,t)\right)\label{d3x3}
\end{equation}
holds. If the displacement satisfies
(\ref{gspecs}) on a time interval $I$ and the local Rossby number $\rho$
does not exceed $1/4$  
then
\begin{equation}
\sup_{t\in I}\|\partial_{a_3}\left (X(a,t) -a\right)\|_{L^{\infty}(da)}\le 9\rho
\label{dzxbound}
\end{equation}
holds on the interval $I$.
\end{prop}
Rearranging (\ref{d3x3}) we have
$$
\left (\frac{\partial X_3}{\partial a_3}(a,t) - 1\right ) + 
\frac{1}{2\Omega\left (1 + \frac{\zeta_3(a)}{2\Omega}\right)} \left (\zeta_1(a)\frac{\partial X_1}{\partial a_3}(a,t) + \zeta_2(a)\frac{\partial X_2}{\partial a_3}(a,t)\right )
$$
\begin{equation}
= \frac{1}{2\Omega\left(1 + \frac{\zeta_3(a)}{2\Omega}\right)}\left (\omega_3(X(a,t),t) - \zeta_3(a)\right)\label{tp}
\end{equation}
Recalling the definition (\ref{lagdis}) of the Lagrangian displacement we have
\begin{thm}
The Lagrangian displacement $\lambda (a,t)$ of a solution of the rotating Euler
equations (\ref{re}) for local Rossby number less than two 
obeys the differential equation
$$
\partial_{a_3}\lambda_3(a,t) + \epsilon_1(a)\partial_{a_3}\lambda_1(a,t)
+ \epsilon_2(a)\partial_{a_3}\lambda_2(a,t) 
$$
\begin{equation}
= \frac{1}{2\Omega\left(1 + \frac{\zeta_3(a)}{2\Omega}\right)}\left (\omega_3(X(a,t),t) - \zeta_3(a)\right)
\label{diseq}  
\end{equation}
where $\epsilon_j(a)$, $j =1, 2$ are given by
\begin{equation}
\epsilon_j(a) = \frac{1}{2\Omega\left (1 + \frac{\zeta_3(a)}{2\Omega}\right)}\zeta_j(a).\label{epsj}
\end{equation}
\end{thm}

It is interesting to note that, if the solution of the
Euler equation is steady, then 
\begin{equation}
\omega(X(a,t),t) = \zeta(a)\label{omegzeta}
\end{equation}
and consequently (\ref{diseq}) becomes 
\begin{equation}
\partial_{a_3}\lambda_3(a,t) + \epsilon_1(a)\partial_{a_3}\lambda_1(a,t)
+ \epsilon_2(a)\partial_{a_3}\lambda_2(a,t) = 0
\label{stdiseq}  
\end{equation}

If the local Rossby number is smaller than two then, if at some label
$a$ the vorticity direction is vertical, then at the same label
the vertical derivative of the vertical Lagrangian displacement vanishes.
The vertical distance between fluid parcels
does not change in time if the parcels were initially on a vertical vortex 
line. 

Let us consider the motion of a rotating incompressible ideal fluid
for a duration $T$, and assume that in this period of time the flow is 
smooth, and so, 
\begin{equation}
\int_0^T\|\nabla u (\cdot,t) \|_{L^{\infty}} <\infty\label{mone}
\end{equation} 
Let us denote
\begin{equation}
g(t) = \|\nabla\ell (\cdot,t)\|_{L^{\infty}}. 
\end{equation}
In view of the equation
\begin{equation}
(\partial_t + u\cdot\nabla)(\nabla\ell) + \nabla u + (\nabla \ell)\nabla u = 0
\label{nablaell}
\end{equation}
obeyed by the gradient of the displacement, we deduce
that
\begin{equation}
g(t) \le \left (e^{\int_{t_0}^t\|\nabla u(\cdot, t)\|_{L^{\infty}(dx)}dt} - 1 \right )\label{gt}
\end{equation}
holds on any interval $I =[t_0, t_0+\tau]$ where $\ell(x,t_0) = 0$. 
Consequently
\begin{equation}
\|(\nabla A_1 \times \nabla A_2)(x,t)\|_{L^{\infty}(dx)} \le e^{2\int_{t_0}^t\|\nabla u(\cdot, t)\|_{L^{\infty}(dx)}dt }
\label{nabels}
\end{equation}
holds on the same time interval, and, in view of (\ref{nabl})
\begin{equation}
\left \| \frac{\partial X}{\partial a_3}(a,t)\right \|_{L^{\infty}(da)} \le e^{2 \int_{t_0}^t\|\nabla u(\cdot, t)\|_{L^{\infty}(dx)} dt}
\label{d3X}
\end{equation}
also holds. 
If wish to ensure that the Lagrangian displacement obeys (\ref{gspecs})
on the interval of time $I = [t_0, t_0 + \tau]$, then,  if (\ref{mone}) holds, 
it is enough to require
\begin{equation}
\int_{t_0}^{t_0 + \tau}\|\nabla u(\cdot, t)\|_{L^{\infty}(dx)}dt \le \log{\frac{5}{4}}.\label{dur}
\end{equation}
Now we are in position to prove that the vertical separation of Lagrangian trajectories is controlled by the local Rossby number:
\begin{thm} Consider a smooth (\ref{mone}) solution of the rotating Euler 
equations (\ref{re}), defined on a time interval $[0, T]$. Consider, at 
time $t_0$, two particles $P =
(a, b, c)$ and $Q = (a,b, c+d)$ separated by a vertical segment of
length $d$. (Vertical refers to the direction of the rotation axis $e_3 = {\widehat z}$).
Consider the Lagrangian evolution of the particles $X(P,t)$, $X(Q,t)$. Then,
the vertical separation between the two particles obeys
\begin{equation}
\left |X_3(P,t)-X_3(Q,t) + d \right | \le \frac{\rho}{2} \left(1 + e^{2\int_{t_0}^t\|\nabla u(\cdot, t)\|_{L^{\infty}(dx)}dt }\right )\label{XX}
\end{equation}
for all $t_0\le t\le T$. 
\end{thm}
For a lower bound we need to restrict the duration of time. We consider
two fluid masses (sets) $\Sigma_j(t)$, $j =1, 2$. We define their
vertical separation as
\begin{equation}
\delta(t) = \inf\{ |z_1-z_2|\,\,; \,\,(x,y,z_j) \in \Sigma_j(t)\}\label{deltat}
\end{equation}
and distance $d(t)$
\begin{equation}
d(t) = \inf\{ |(x_1,y_1,z_1) - (x_2,y_2,z_2)\,\, ;\,\, (x_j,y_j,z_j)\in \Sigma_j(t)\}\label{dist}
\end{equation}
Obviously, $d(t)\le\delta(t)$.  
\begin{thm}
Let $\Sigma_j(t)$ be two sets carried by the flow of a smooth (\ref{mone}) 
solution of the rotating Euler equations (\ref{re}). Assume
that at some time $t_0$, the distance between the two sets is positive,
 $d(t_0)>0$. Consider a time interval $I=[t_0, t_0+\tau]$ in which
(\ref{gspecs}) holds. Then
\begin{equation}
\delta(t) \ge (1 + 14\rho)^{-1}d(t_0)\label{lsepb}
\end{equation}
holds for $t\in [t_0, t_0 +\tau]$
\end{thm}
 
Indeed, if two points are in the respective sets, $(x,y,z_j)\in \Sigma_{j}(t)$,
then, considering the back-to-labels map starting from time
$t_0$, $\ell (x,t_0) =0$ we have by assumption $A(x,y,z_j, t)\in \Sigma_j(t_0)$, and, in view of (\ref{dzlbound}) we have
$$
d(t_0) \le |A(x,y,z_1,t)-A(x,y,z_2,t)| = |z_1-z_2 + \ell(x,y,z_1,t)-
\ell(x,y,z_2,t)| \le 
$$
$$
|z_1-z_2|(1 + 14 \rho).
$$
Taking the infimum produces (\ref{lsepb}).
The bounds above hold in great generality; they become most effective when
$\rho\to 0$. It is quite easy to exhibit at least some situations in which
this limit is well defined mathematically.
\begin{prop}
Consider the rotating Euler system (\ref{re}) with spatially periodic
boundary conditions,
$u(x+L_1,y+L_2, z+L_3,t) = u(x,y,z,t)$. Assume that the initial data
$u(x,y,z,0) = u_0(x,y,z)$ belongs to the space $H^s$ of divergence-free 
square integrable functions having $s>\frac{5}{2}$ square integrable derivatives. Then, there exists a time $T$ and a constant $M_1$, depending only on 
the norm of the initial data in $H^s$, such that, for any $\Omega>0$, the 
solution of (\ref{re}) with initial data $u_0$ exists on the time 
interval $[0,T]$ and obeys
\begin{equation}
\sup_{0\le t\le T}\|\nabla u(\cdot, t)\|_{L^{\infty}(dx)}\le M_1.\label{m1}
\end{equation}
\end{prop}
By restricing the duration of time further (see \ref{dur}),
\begin{equation}
M_1T \le \log{\frac{5}{4}}\label{mt}
\end{equation}
we guarantee (\ref{gspecs}) on the time interval $[0,T]$. Fixing
thus $u_0$ and $T$, we may let $\Omega\to\infty$ and consequently
$\rho \le\frac{M_1}{\Omega}\to 0$. Assuming that
\begin{equation}
\Omega \ge 4 M_1\label{omegab}
\end{equation}
we have the conditions for the bounds
(\ref{dzlbound}, \ref{dzxbound}) to hold.
\begin{thm}
Let $u_0\in H^s$, $s>\frac{5}{2}$ and $T>0$ be fixed as above 
satisfying (\ref{m1}, \ref{mt}). Then, for
each $\Omega$ satisfying (\ref{omegab}), consider the inverse and
direct Lagrangian displacements 
$$
\ell (x,y,z,t) = A(x,y,z,t) -(x,y,z)
$$ 
and 
$$
\lambda (a_1,a_2,a_3,t) = X(a_1,a_2,a_3,t) -(a_1,a_2,a_3).
$$
They obey
\begin{equation}
\left \|\partial_z \ell(\cdot, t)\right \|_{L^{\infty}(dx)} \le 14\rho
\label{dzl}
\end{equation}
and
\begin{equation}
\left \|\partial_{a_3} \lambda (\cdot, t)\right \|_{L^{\infty}(da)}                                                       \le 9\rho
\label{dzlambda} 
\end{equation}
with $\rho = \Omega^{-1}\sup_{0\le t \le T}\|\omega (\cdot, t)\|_{L^{\infty}(dx)}\le \Omega^{-1}M_1$.

Let $\Omega_j\to\infty$ be an arbitrary sequence and
let $X^{j}(a_1,a_2,a_3,t)$ denote the Lagrangian paths
associated to $\Omega_j$. Then, there exists a subsequence
(denoted for convenience by the same letter $j$)
an invertible map $X(a_1,a_2,a_3, t)$, and a periodic function
of two variables $\lambda(a_1,a_2,t)$  
such that
$$
\lim_{j\to\infty} X^{j}(a_1,a_2,a_3,t) = X(a_1,a_2,a_3, t)
$$
holds uniformly in $a,t$
and
$$
X(a_1,a_2,a_3) = (a_1,a_2,a_3) +\lambda (a_1, a_2, t)
$$

\end{thm}

The proof of the second part follows easily from the Arzeli-Ascola theorem and
inequalities of the form
$$
\left |X^{(j)}(a_1, a_2, b_3,t) - X^{(j)}(a_1, a_2, c_3, t) + c_3-b_3\right |
\le 9\rho_j|b_3-c_3|
$$ 

\noindent{\bf{Acknowledgment}} Research supported in part by NSF DMS-0202531
and by the ASCI Flash Center at the University of Chicago under DOE
contract B341495.

\end{document}